\documentclass[10pt, twoside]{article}

\usepackage{amsmath}
\usepackage{amsfonts}
\usepackage{amssymb}
\usepackage{eucal}

\newtheorem{tw}{Theorem}

\newtheorem{col}{Corollary}

\newtheorem{lem}{Lemma}

\newenvironment{pr}{\noindent{\bf Proof. }}{\vskip5mm}

\newcounter{ex}
\newenvironment{ex}{\stepcounter{ex}\vskip3mm\noindent \bf Example \arabic{ex}.\rm\hskip2mm}{\vskip 0.3cm}

\begin{document}

\thispagestyle{empty}
\begin{center}
\textbf{PARA-CR STRUCTURES \\ON ALMOST PARACONTACT METRIC MANIFOLDS}
\end{center}
\vskip 1cm
\begin{center}
\footnotesize{
JOANNA WE{\L}YCZKO\\[4pt]
{\it Institute of Mathematics and Computer Science,\\
Wroc{\l}aw University of Technology,\\
Wybrze\.{z}e Wyspia\'nskiego 27, 50-370 Wroc{\l}aw, Poland}\\
{\it E-mail}: {\it Joanna.Welyczko}{\char64}{\it pwr.wroc.pl}}
\end{center}

\bigskip
\begin{minipage}{4.5 in}
\baselineskip 10pt
\footnotesize{ 
 Almost paracontact metric manifolds are the famous examples of almost para-CR manifolds. We find necessary and sufficient conditions for such manifolds for be para-CR. Next we examine these conditions in certain subclasses of almost paracontact metric manifolds. Especially, it is shown that the normal almost paracontact metric manifolds are para-CR. We establish necessary and sufficient conditions for paracontact metric manifolds as well as for almost paracosymplectic manifolds to be para-CR. We find also basic curvature identities for para-CR paracontact metric manifolds and study their consequences. Among others, we prove that any para-CR paracontact metric manifold of constant sectional curvature and of dimension greather than tree must be para-Sasakian and its curvature equal to minus one. The last assertion do not hold in dimension tree. Moreover, we show that a conformally flat para-Sasakian manifold  is of constant sectional curvature equal to minus one. 
New classes of examples of para-CR manifolds are established.
 \\[+6pt]
}

\end{minipage}
\baselineskip 13pt

{\section*{\normalsize{1. Preliminaries}}}

\noindent
Let $M$ be an almost paracontact manifold and $(\varphi,\xi,\eta)$ its almost paracontact  structure (e.g. \cite{E,SZ}). This means that $M$ is an $(2n+1)$-dimensional differentiable manifold and $\varphi, \xi, \eta$ are tensor fields on $M$ of type $(1,1)$, $(1,0)$, $(0, 1)$, respectively, such that 
$$
  \varphi^2 X=X-\eta(X)\xi,\quad \eta(\xi)=1,\quad \varphi\xi=0,\quad \eta\circ\varphi=0. 
$$
Moreover the tensor field $\varphi$ induces an almost paracomplex structure on the paracontact distribution $\mathcal{D}=\mathop{\rm Ker}\eta$, i.e. the eigendistributions $\mathcal{D}^{\pm}$ corresponding to the eigenvalues $\pm1$ of $\varphi$ are both $n$-dimensional.
A pseudo-Riemannian metric $g$ on $M$ satisfying the condition
$$
  g(\varphi X,\varphi Y)=-g(X,Y)+\eta(X)\eta(Y)
$$
is said to be compatible with the structure $(\varphi, \xi, \eta)$. If $g$ is such a metric, then the quadruplet $(\varphi,\xi,\eta,g)$ is called an almost paracontact  metric structure and $M$ an almost paracontact metric  manifold. For such a manifold, we additionally have $\eta(X)=g(X,\xi)$, and we define the (skew-symmetric) fundamental 2-form $\varPhi$ by $\varPhi(X, Y)=g(X,\varphi Y)$. 

In the above and in the sequel, $W,X,Y,Z,...$ indicate arbitrary vector fields on the considered manifold if it is not otherwise stated. 

An almost paracontact metric manifold  is called to be \\
\rm{(a)} normal if (\cite{K,KW})
$$
  N(X,Y)-2\,d\eta(X,Y)\xi=0,
$$
\hskip 0.6cm where $N$ is the Nijenhuis torsion tensor of $\varphi$,
$$
  N(X,Y)=\varphi^2[X,Y]+[\varphi X,\varphi Y]-\varphi([\varphi X,Y]+[X,\varphi Y]);
$$
\rm{(b)} paracontact metric if $\varPhi=d\eta$ (\cite{E,SZ});\\
\rm{(c)} para-Sasakian if is normal and paracontact metric;\\
\rm{(d)} almost para-cosympletic if the forms $\eta$ and $\varPhi$ are closed, that is, $d\eta=0$ and 
$ d\varPhi=0$ (\cite{D}).

\vspace*{5mm}
\section*{\normalsize 2. Para-CR Manifolds}

\noindent
Almost paracontact metric manifolds can be interpreted as almost para-CR manifolds. Following \cite{AMT1,AMT2,CM,IVZ,SZ}, we will resume the demanded details. It is also important to mention certain analytic and geometric studies on generalizations of para-CR structures which have occured in the papers \cite{HN,Kush,NS}.

Let $M$ be an almost paracontact metric manifold and $(\varphi,\xi,\eta,g)$ its almost paracontact metric structure. 
Then $\dim\mathcal{D}=2n$ and $\varphi$ (precisely, $\varphi|_\mathcal{D})$ is a field of endomorphisms of $\mathcal{D}$ such that $\varphi^2=\mathop{\rm Id}$, and the eigendistributions $\mathcal{D}^{\pm}$ corresponding to the eigenvalues $\pm1$ of $\varphi$ are both $n$-dimensional. Thus, the pair $(\mathcal{D},\varphi)$ becomes an almost para-CR structure on $M$. 

We say that $(\mathcal{D},\varphi)$ is a para-CR structure if it is formally integrable, that is, the following two conditions are satisfied 
\begin{eqnarray}
\label{s0}
  & [\varphi X,Y]+[X,\varphi Y]\in\mathcal{D}, & \\
\label{s1}
  &S(X,Y)= [X,Y]+[\varphi X,\varphi Y]-\varphi([X,\varphi Y]+[\varphi X,Y])=0 &
\end{eqnarray}
for all $X,Y\in\mathcal{D}$. Equivalently, the formal integrability means that the eigendistributions $\mathcal{D}^{+}$ and $\mathcal{D}^{-}$ are involutive, that is, 
\begin{equation}
\label{plusminus}
  [\mathcal{D}^{+},\mathcal{D}^{+}]\subset \mathcal{D}^{+},\quad
  [\mathcal{D}^{-},\mathcal{D}^{-}]\subset \mathcal{D}^{-}.
\end{equation}
$M$ will be called to be a para-CR manifold if $(\mathcal{D},\varphi)$ is a para-CR structure.

In the sequel, we need certain new shapes of the condition (\ref{s0}). 

\begin{lem}
\label{auxlem}
Any of the following conditions is equivalent to (\ref{s0})
\begin{eqnarray}
\label{news00}
  & d\eta(X,\varphi Y) - d\eta(Y,\varphi X) = 0, & \\
\label{news01}
  & (\nabla_X\eta)(\varphi Y) + (\nabla_{\varphi X}\eta)(Y)
     = (\nabla_Y\eta)(\varphi X) + (\nabla_{\varphi Y}\eta)(X) &
\end{eqnarray}
for any $X,Y\in\mathcal D$.
\end{lem}

\begin{pr}
Indeed, note that (\ref{s0}) is fulfilled if and only if $\eta([\varphi X,Y]+[X,\varphi Y])=0$. Now, since $d\eta(U,V)=-(1/2)\eta([U,V])$ for any $U,V\in\mathcal D$, we claim easily that the condition (\ref{news00}) is equivalent to (\ref{s0}). On the other hand, since in general
$$
  d\eta(U,V)=\frac12((\nabla_U\eta)(V) - (\nabla_V\eta)(U)),
$$
the condition (\ref{news00}) can be equivalently written as (\ref{news01}).
\hfill $\Box$
\end{pr}

When we define a $(0,2)$-tensor field $L$ on $\mathcal{D}$ by 
\begin{equation}
\label{levi}
  L(X,Y)=-d\eta(X,\varphi Y),\quad X,Y\in\mathcal{D},
\end{equation}
then the condition (\ref{news00}) can be interpreted as the symmetry of $L$. By an analogy to the theory of CR structures (cf. e.g. \cite{DT}), the tensor $L$ defined by (\ref{levi}) can be called the Levi form corresponding to the para-CR structure. 

The following theorem gives a necessary and sufficient condition for $M$ to be para-CR.

\begin{tw}
\label{tw1}
An almost paracontact metric manifold $M$ is a para-CR manifold if and only if the following condition
\begin{equation}
\label{paracr}
  (\nabla_X\varphi)Y+(\nabla_{\varphi X}\varphi)\varphi Y 
     =\null-((\nabla_Y\eta)(\varphi X)+(\nabla_{\varphi Y}\eta)(X))\xi 
\end{equation}
is satisfied for any $X,Y\in \mathcal{D}$, with $\nabla$ being the Levi-Civita connection of $M$.
\end{tw}

\begin{pr}
Before we start with the proof, it will be useful to find the following expression 
\begin{eqnarray}
\label{s2}
  \varphi S(X,Y) 
  &=& \varphi([X,Y]+[\varphi X,\varphi Y])-[X,\varphi Y]-[\varphi X,Y]\nonumber \\
  & & \null+\eta([X,\varphi Y]+[\varphi X,Y]\xi \nonumber\\
  &=& \varphi(\nabla_{\varphi X}\varphi)Y-\varphi(\nabla_{\varphi Y}\varphi)X
      -(\nabla_X\varphi)Y+(\nabla_Y\varphi)X \nonumber\\
  & & \null+\eta([X,\varphi Y]+[\varphi X,Y])\xi \nonumber\\
  &=& \null-(\nabla_{\varphi X}\varphi)\varphi Y
      +(\nabla_{\varphi Y}\varphi)\varphi X-(\nabla_X\varphi)Y+(\nabla_Y\varphi)X \nonumber\\
  & & \null+\eta([X,\varphi Y]+[\varphi X,Y])\xi 
\end{eqnarray}
for any $X,Y\in\mathcal D$.

Let us assume that (\ref{s0}) and (\ref{s1}) are satisfied. Define an auxiliary $(0,3)$-tensor field $A$ on $\mathcal{D}$ by 
\begin{equation}
\label{dfA}
  A(X,Y,Z)=g((\nabla_{\varphi_X}\varphi)\varphi Y+(\nabla_X\varphi)Y,Z)
\end{equation}
for any $X,Y,Z\in\mathcal D$. Applying (\ref{s2}) with $S=0$, we claim that 
\begin{equation}
\label{sym}
  A(X,Y,Z) = A(Y,X,Z).
\end{equation}
Moreover, by a simply calculation we show that
\begin{eqnarray}
\label{antysym}
  && A(X,Y,Z) + A(X,Z,Y)\nonumber \\
  &&\quad=g((\nabla_{\varphi X}\varphi)\varphi Y+(\nabla_{X}\varphi)Y,Z)
      +g((\nabla_{\varphi X}\varphi)\varphi Z+(\nabla_{X}\varphi)Z,Y)\nonumber\\
  &&\quad=\null-g((\nabla_{\varphi X}\varphi)Z,\varphi Y)
               -g((\nabla_{\varphi X}\varphi)Y,\varphi Z)=0.
\end{eqnarray}
Using (\ref{sym}) and (\ref{antysym}), we can compute
\begin{eqnarray*}
 A(X,Y,Z) &=& -A(X,Z,Y) = -A(Z,X,Y) = A(Z,Y,X) \\
          &=& A(Y,Z,X) = -A(Y,X,Z) = -A(X,Y,Z).
\end{eqnarray*}
Hence, $A(X,Y,Z)=0$ for any $X,Y,Z\in\mathcal{D}$. By the definition (\ref{dfA}), this implies the following
\begin{equation}
\label{lambda}
  (\nabla_{\varphi X}\varphi)\varphi Y+(\nabla_X\varphi)Y=\lambda(X,Y)\xi\quad {\rm for}\; X,Y\in \mathcal{D},
\end{equation}
for a certain $(0,2)$-tensor field $\lambda$ on $\mathcal D$. The projection (\ref{lambda}) onto $\xi$ leads to
\begin{eqnarray}
\label{lamb}
\lambda(X,Y)
  &=&g((\nabla_{\varphi X}\varphi)\varphi Y+(\nabla_X\varphi)Y,\xi)\nonumber\\
  &=&\null-g((\nabla_{\varphi X}\varphi)\xi,\varphi Y)-g((\nabla_X\varphi)\xi,Y)\nonumber\\
  &=&g(\varphi\nabla_{\varphi X}\xi,\varphi Y)-g((\nabla_X\varphi)\xi,Y)\nonumber\\
  &=&\null-g(\nabla_{\varphi X}\xi,Y)+g(\varphi\nabla_X\xi,Y)\nonumber\\
  &=&g(\varphi\nabla_X\xi-\nabla_{\varphi X}\xi,Y)\nonumber\\
  &=&\null-(\nabla_X\eta)\varphi Y-(\nabla_{\varphi X}\eta)Y.
\end{eqnarray} 
By applying (\ref{news01}) (which is equivalent to (\ref{s0})) into (\ref{lamb}), we see that $\lambda(X,Y)$ can be written as
$$
  \lambda(X,Y) = \null-(\nabla_Y\eta)\varphi X-(\nabla_{\varphi Y}\eta)X,
$$
which substituted into (\ref{lambda}), gives (\ref{paracr}).

Conversely, assume that (\ref{paracr}) is fulfilled. Projecting (\ref{paracr}) onto the vector field $\xi$, we obtain (\ref{news01}), and consequently (\ref{s0}) too. But using (\ref{news01}) and (\ref{paracr}), we have 
$$
  (\nabla_{\varphi X}\varphi)\varphi Y+(\nabla_X\varphi)Y
   = (\nabla_{\varphi Y}\varphi)\varphi X+(\nabla_Y\varphi)X.
$$
Now, using the above and (\ref{s0}) into (\ref{s2}), we have $\varphi S(X,Y)=0$. Hence, $S(X,Y)=0$. Thus, we get (\ref{s1}), completing the proof.
\hfill $\Box$
\end{pr}

\begin{tw}
Any 3-dimensional almost paracontact metric manifold is a para-CR manifold.
\end{tw}

\begin{pr}
Recall that it is proved by the author in \cite{JW} that for any 3-dimensional almost paracontact metric manifold, it holds 
$$
  (\nabla_X \varphi)Y=g(\varphi\nabla_X \xi,Y)\xi-\eta(Y)\varphi\nabla_X \xi
$$
for any $X,Y\in\mathfrak{X}(M)$. Now, using the above formula, we get for such a manifold,
\begin{eqnarray*}
(\nabla_X\varphi)Y+(\nabla_{\varphi X}\varphi)\varphi Y
  &=&g(\varphi\nabla_X \xi,Y)\xi-g(\nabla_{\varphi X} \xi,Y)\xi \\
  &=&\null-((\nabla_X\eta)(\varphi Y)+(\nabla_{\varphi X}\eta) (Y))\xi
\end{eqnarray*}
for any $X,Y\in \mathcal{D}$. In view of the above and Theorem \ref{tw1}, the condition (\ref{paracr}) reduces to (\ref{news01}), or equivalently, to (\ref{news00}). But since $\dim\mathcal D=2$, the verification of (\ref{news01}) is easy. In fact it is sufficient to take $X=E_1$ and $Y=E_2$, where $E_1$, $E_2$ form a local basis for $\mathcal D$ such that $\varphi E_1=-E_1$ and $\varphi E_2=E2$.
\hfill $\Box$
\end{pr}

\section*{\normalsize 3.  Normal Almost Paracontact Metric Manifolds}

\noindent
We start with recalling the theorem proved by the author in \cite{JW}: An almost paracontact metric manifold is normal if and only if
\begin{equation}
\label{normal}
  \varphi(\nabla_X\varphi)Y-(\nabla_{\varphi X}\varphi)Y+(\nabla_X\eta)(Y)\xi=0
\end{equation}
for any $X,Y\in\mathfrak{X}(M)$. 

Assume that $M$ is a normal almost paracontact metric manifold. For such a manifold, we additionally have the following relations 
\begin{eqnarray}
\label{wlasn1}
  & \nabla_{\xi}{\xi} = 0,\quad \nabla_{\xi}\eta=0, & \\
\label{wlasn2}
  & \nabla_{\varphi X}\xi = \varphi\nabla_X\xi, & \\
\label{wlasn3}
  & \nabla_{\xi}\varphi = 0. &
\end{eqnarray}
Indeed, by putting $X=Y=\xi$ into the formula (\ref{normal}), we obtain $0 = \varphi(\nabla_{\xi}\varphi)\xi = -\varphi^2\nabla_{\xi}\xi$. Hence, $\nabla_{\xi}\xi=0$ and $\nabla_{\xi}\eta=0$. Putting $Y=\xi$ into (\ref{normal}), we get 
$$
  0 = \varphi(\nabla_X\varphi)\xi - (\nabla_{\varphi X}\varphi)\xi
    = -\nabla_X\xi + \varphi\nabla_{\varphi X}\xi.
$$
Hence, (\ref{wlasn2}) follows. Putting $X=\xi$ into (\ref{normal}) and using (\ref{wlasn1}), we obtain $\varphi(\nabla_{\xi}\varphi)Y=0$, which implies (\ref{wlasn3}). 

The following theorem is the main result of this section:

\begin{tw}
\label{N}
Any normal almost paracontact metric manifold is a para-CR manifold.
\end{tw}

\begin{pr}
Let $M$ be a normal almost paracontact metric manifold. Let us suppose that $X,Y\in\mathcal{D}$. From (\ref{normal}), we deduce
\begin{equation}
\label{aux1}
  \varphi(\nabla_X\varphi)\varphi Y-(\nabla_{\varphi X}\varphi)\varphi Y+(\nabla_X\eta)(\varphi Y)\xi=0.
\end{equation}
On the other hand, it is easy to see that
$$
  \varphi(\nabla_X\varphi)\varphi Y = \null-(\nabla_X\varphi)Y - (\nabla_X\eta)(\varphi Y)\xi,
$$
which together with (\ref{aux1}) leads to
\begin{equation}
\label{aux2}
  (\nabla_X\varphi)Y + (\nabla_{\varphi X}\varphi)\varphi Y=0.
\end{equation}
Moreover, we can show that 
\begin{equation}
\label{aux3}
  (\nabla_{\varphi Y}\eta)(X)+(\nabla_Y\eta)(\varphi X) = 0.
\end{equation}
In fact, using (\ref{wlasn2}), the left hand side of (\ref{aux3}) can be transformed in the following way
$$
  g(\nabla_{\varphi Y}\xi,X)+g(\nabla_Y\xi,\varphi X) 
  = g(\varphi\nabla_Y\xi,X)+g(\nabla_Y\xi,\varphi X)=0.
$$
Finally, having (\ref{aux2}) and (\ref{aux3}), we claim that the relation (\ref{paracr}) is obviously fulfilled. In view of Theorem \ref{tw1}, this completes the proof. 
\hfill $\Box$
\end{pr}

\section*{\normalsize 4.  Paracontact Metric Manifolds}

\noindent
We recall certain facts about paracontact metric manifolds from the paper \cite{SZ}.

Let $M$ be a paracontact metric manifold. Define $h=\frac12\mathcal{L}_\xi\varphi$, where $\mathcal{L}$ indicates the Lie derivation operator. Then, 
\begin{eqnarray}
\label{h1}
  & g(hX,Y)=g(hY,X)\quad\mbox{($h$ is a symmetric operator)}, & \\
\label{h2}
  & \varphi h+h\varphi=0,\quad \mathop{\rm Tr} h=0,\quad h\xi=0,\quad \eta\circ h=0. 
\end{eqnarray}

It is very important that on a paracontact metric manifold $M$, the following relations hold
\begin{eqnarray}
\label{h}
  \nabla_X\xi &=&\null-\varphi X+\varphi hX, \\
\label{lemat}
  \quad(\nabla_{\varphi X}\varphi)\varphi Y-(\nabla_X\varphi)Y
  &=& 2g(X,Y)\xi-\eta(Y)(X-hX+\eta(X)\xi). 
\end{eqnarray}

A paracontact metric manifold is a para-Sasakian one if and only if 
\begin{equation}
\label{sas}
(\nabla_X\varphi)Y=-g(X,Y)\xi+\eta(Y)X, 
\end{equation}
and on a para-Sasakian manifolds, it always holds
\begin{equation}
\label{pScons}
  \nabla_X\xi = -\varphi X, \quad h = 0.
\end{equation}

\begin{tw}
\label{K}
A paracontact metric manifold is a para-CR manifold if and only if 
\begin{equation}
\label{wzor1}
  (\nabla_X\varphi)Y=g(\varphi\nabla_X\xi,Y)\xi-\eta(Y)\varphi\nabla_X\xi.
\end{equation}
\end{tw}

In \cite[Th.\ 4.10]{SZ}, S. Zamkovoy  proved that the paracontact metric manifold is para-CR manifold if and only if
\begin{equation}
\label{wzorzamk}
  (\nabla_X\varphi)Y = -g(X-hX,Y)\xi + \eta(Y)(X-hX).
\end{equation}
In view of (\ref{h}), the condition  (\ref{wzorzamk}) is equivalent to (\ref{wzor1}).
 
However, using our convention we shall give a short prove of this fact.

\begin{pr}
Let $M$ be a paracontact manifold. For such a manifold, using (\ref{h}) and (\ref{h2}), we obtain for $X,Y\in\mathcal{D}$, 
\begin{equation}
\label{cont1}
  (\nabla_Y\eta)(\varphi X) + (\nabla_{\varphi Y}\eta)(X) = -2g(hY,X)\xi=-2d(hX,Y).
\end{equation}
Moreover, from (\ref{lemat}), for $X,Y\in\mathcal{D}$, we get 
\begin{equation}
\label{cont2}
  (\nabla_{\varphi X}\varphi)\varphi Y = (\nabla_X\varphi)Y + 2g(X,Y)\xi.
\end{equation}
By applying (\ref{cont1}) and (\ref{cont2}) into (\ref{paracr}), we claim from Theorem \ref{tw1} that the paracontact metric manifold $M$ is a para-CR manifold if and only if 
\begin{equation}
\label{contparacr}
  (\nabla_X\varphi)Y = -g(X-hX,Y)\xi
\end{equation}
for any $X,Y\in \mathcal{D}$. 

It is clear that (\ref{wzorzamk}) implies (\ref{contparacr}). Conversely, take arbitrary vector fields $X,Y\in\mathfrak{X}(M)$. Then $X-\eta(X)\xi\in\mathcal{D}$ and $Y-\eta(Y)\xi\in\mathcal{D}$. Put $X-\eta(X)\xi$ instead of $X$ and $Y-\eta(Y)\xi$ instead of $Y$ into the formula (\ref{contparacr}). After some calculations, in which the above listed properties of $\nabla\varphi$ and $h$ should be used, one obtains (\ref{wzorzamk}). This completes the proof. 
\hfill $\Box$
\end{pr}

A para-Sasakian manifold is obviously para-CR paracontact metric.

\begin{col}
\label{wnsas}
A para-CR paracontact metric manifold is para-Sasakian if and only if $h=0$.
\end{col}

\begin{pr}
Comparing (\ref{sas}) with (\ref{wzorzamk}), we get
$$
g(hX,Y)\xi-\eta(Y)hX=0,
$$
which is equivalent to $h=0$.
\hfill $\Box$
\end{pr}

\begin{tw}
For a para-CR paracontact metric manifold, we have the following curvature identity
\begin{eqnarray}
\label{k1}
  && \quad (R(W,X)\varphi)Y = R(W,X)\varphi Y - \varphi R(W,X)Y\nonumber \\
  && \qquad =g((\nabla_Wh)X-(\nabla_Xh)W,Y)\xi\nonumber\\
  && \quad\qquad\null+g(hX-X,Y)\varphi(hW-W)-g(hW-W,Y)\varphi(hX-X)\nonumber\\
  && \quad\qquad\null-g(\varphi(hW-W),Y)(hX-X)+g(\varphi(hX-X),Y)(hW-W)\nonumber\\
  && \quad\qquad\null-\eta(W)((\nabla_Wh)X-(\nabla_Xh)W).
\end{eqnarray}
\end{tw}

\begin{pr}
Using (\ref{wzor1}) and (\ref{h}), we calculate
\begin{eqnarray*}
  (\nabla^2_{WX}\varphi)Y
  &=&\nabla_W((\nabla_X\varphi)Y)-(\nabla_{\nabla_WX}\varphi)Y-(\nabla_X\varphi)\nabla_WY\\
  &=& \nabla_W(g(hX-X,Y)\xi-\eta(Y)(hX-X))\\
  & & \null-g(h(\nabla_WX)-\nabla_WX,Y)\xi+\eta(Y)(h\nabla_WX-\nabla_WX)\\
  & & \null-g(hX-X,\nabla_WY)\xi+\eta(\nabla_WY)(hX-X)\\
  &=& g((\nabla_Wh)X,Y)+g(hX-X,Y)\varphi(hW-W)\\
  & & \null-g(Y,\varphi(hW-W))(hX-X)-\eta(Y)(\nabla_Wh)X.
\end{eqnarray*}
We obtain the final assertion by applying the above formula into the below identity
$$
  (R(W,X)\varphi)Y=(\nabla^2_{WX}\varphi)Y-(\nabla^2_{XW}\varphi)Y.
$$
\hfill $\Box$
\end{pr}

\begin{col}
For a paracontact metric para-CR manifold, we have
\begin{eqnarray}
\label{k2}
  \quad g(R(W,X)\varphi Y,\xi)
     =g(\nabla_Wh)X-(\nabla_Xh)W,Y)-2\eta(Y)g(\varphi h^2W,X) &&\nonumber \\
     \null+\eta(X)g(\varphi(hW-W),Y)-\eta(W)g(\varphi(hX-X),Y). && 
\end{eqnarray}
\end{col}

\begin{pr}
The projection of (\ref{k1}) onto $\xi$ leads to
\begin{eqnarray}
\label{xxx}
  && \quad g(R(W,X)\varphi Y,\xi)=g((\nabla_Wh)X-(\nabla_Xh)W,Y) \nonumber\\
  && \qquad\qquad\null+\eta(X)g(\varphi(hW-W),Y)-\eta(W)g(\varphi(hX-X),Y)\nonumber\\
  && \qquad\qquad\null-\eta(Y)g((\nabla_Wh)X-(\nabla_Xh)W,\xi).
\end{eqnarray}
On the other hand, using (\ref{h}), we find 
\begin{eqnarray*}
  g((\nabla_Wh)X-(\nabla_Xh)W,\xi)&=&-g(h\nabla_W\xi,X)+g(h\nabla_X\xi,W)\\
  &=&2g(\varphi h^2W,X),
\end{eqnarray*}
which turns (\ref{xxx}) into (\ref{k2}).
\hfill $\Box$
\end{pr}

In \cite[Th.\ 3.12]{SZ}, S. Zamkovoy proved that if a paracontact metric manifold $M$ is of constant sectional curvature $k$ and of dimension $2n+1\geqslant5$, then  $k=-1$ and $|h|^2=0$. From the proof of his theorem it even more follows that $h^2=0$. Unfortunetly, $h^2=0$ does not imply $h=0$, and therefore the manifold does not have to be a para-Sasakian one (see Example \ref{p1}). We will strengthen the assertion of Zamokovoy's theorem making additional assumption that the paracontact metric manifold is para-CR.  

\begin{tw}
A para-CR paracontact metric manifold of constant sectional curvature $k$ and of dimension $2n+1\geqslant5$ is para-Sasakian and $k=-1$.
\end{tw}

\begin{pr}
Since the manifold is of constant sectional curvature $k$, we have
\begin{equation}
\label{yyy}
  R(W,X)\varphi Y = k(g(X,\varphi Y)W-g(W,\varphi Y)X),
\end{equation}
which applied into (\ref{k2}) gives
\begin{eqnarray*}
  \quad g((\nabla_Wh)X-(\nabla_Xh)W,Y)=(k+1)(\eta(W)g(X,\varphi Y)-\eta(X)g(W,\varphi Y)) && \\
  \null+\eta(W)g(\varphi hX,Y)-\eta(X)g(\varphi hW,Y)+2\eta(Y)g(\varphi h^2W,X). &&
\end{eqnarray*}
Hence
\begin{eqnarray}
\label{k3}
  \quad (\nabla_Wh)X-(\nabla_Xh)W=-(k+1)(\eta(W)\varphi X-\eta(X)\varphi W) &&\nonumber \\
      \null+\eta(W)\varphi hX-\eta(X)\varphi hW +2g(\varphi h^2W,X)\xi. &&
\end{eqnarray}
Now, we use (\ref{yyy}) and (\ref{k3}) in (\ref{k1}) and  get
\begin{eqnarray}
\label{k4}
  && \ k(g(X,\varphi Y)W-g(W,\varphi Y)X-g(X,Y)\varphi W+g(W,Y)\varphi X) \nonumber\\
  && \quad= \eta(X)(g((k+1)\varphi W-\varphi hW,Y)\xi-\eta(Y)((k+1)\varphi W-\varphi hW)) \nonumber\\
  && \qquad\null-\eta(W)(g((k+1)\varphi X-\varphi hX,Y)\xi-\eta(Y)((k+1)\varphi X-\varphi hX)) \nonumber\\
  && \qquad\null+g(hX-X,Y)\varphi(hW-W)-g(hW-W,Y)\varphi(hX-X) \nonumber\\
  && \qquad\null+g(\varphi(hX-X),Y)(hW-W)-g(\varphi(hW-W),Y)(hX-X). 
\end{eqnarray}
We calculate the trace of (\ref{k4}) with respect to the arguments $Y$, $W$ and the metric $g$. Then we obtain 
\begin{eqnarray}
\label{k5}
  (2n-2)(k+1)\varphi X=(2n-2)\varphi hX-\mathop{\rm Tr}(\varphi h)(hX-\varphi^2X).
\end{eqnarray}
The left hand side of formula (\ref{k5}) is an antisymmetric linear operator, whereas in view of (\ref{h1}) and (\ref{h2}) the right hand side of (\ref{k5}) is a symmetric operator. Hence, we infer that 
\begin{eqnarray}
\label{k6}
  & (2n-2)(k+1)\varphi X=0, & \\
\label{k8}
  & (2n-2)\varphi hX-\mathop{Tr}(\varphi h)(hX-\varphi^2X)=0. &
\end{eqnarray}
Because of $n>1$, from (\ref{k6}), we have $k=-1$. Moreover, calculating the trace of (\ref{k8}) and using (\ref{h2}), we get $\mathop{\rm Tr}(\varphi h)=0$. Therefore, (\ref{k8}) takes the form $\varphi hX=0$. Hence $h=0$ and by Collorary  \ref{wnsas}, the manifold is para-Sasakian. 
\hfill $\Box$ 
\end{pr}

As it follows from the following example, the assertion of the above theorem does not hold in dimension 3. 

\begin{ex}
Define an almost paracontact metric structure $(\varphi,\xi,\eta,g)$ on $\mathbb R^3$ by assuming 
\begin{eqnarray*}
  & \varphi\dfrac{\partial}{\partial x} = \cosh(2z)\dfrac{\partial}{\partial z},\quad 
    \varphi\dfrac{\partial}{\partial y} = \sinh(2z)\dfrac{\partial}{\partial z}, &\\ 
  & \varphi\dfrac{\partial}{\partial z} = \cosh(2z)\dfrac{\partial}{\partial x}
            -\sinh(2z)\dfrac{\partial}{\partial y}, &\\
  & \xi = \null-\sinh(2z)\dfrac{\partial}{\partial x}+\cosh(2z)\dfrac{\partial}{\partial y},\quad 
    \eta = \sinh(2z)dx + \cosh(2z)dy, &\\
  & g = \null-dx\otimes dx+dy\otimes dy+dz\otimes dz, &
\end{eqnarray*}
$(x,y,z)$ being the Cartesian coordinates in $\mathbb R^3$. In fact, this structure is flat, tree-dimensional, para-CR and paracontact metric. Moreover, we can compute for the tensor field $h$,
\begin{eqnarray*}
  h\frac{\partial}{\partial z} 
  &=& \frac12 \big(\mathcal L_{\xi}\varphi\big)\frac{\partial}{\partial z}
      = \frac12 \Big(\mathcal L_{\xi}\varphi\frac{\partial}{\partial z} 
        - \varphi\mathcal L_{\xi}\frac{\partial}{\partial z}\Big) \\
  &=& \frac12 \Big(\Big[\xi,\varphi\frac{\partial}{\partial z}\Big] 
      - \varphi\Big[\xi,\frac{\partial}{\partial z}\Big]\Big)
      = -\frac{\partial}{\partial z},
\end{eqnarray*}
which shows that $h\not=0$, and therefore the structure is not para-Sasakian. 
\hfill$\Box$
\end{ex}

A para-Sasakian structure of any dimension $2n+1\geqslant3$ and of constant sectional curvature equal to 
minus one can be constructed in the following way.

\begin{ex}
\label{sfera}
Let $(J,G)$ be the standard flat para-K\"ahler structure on $\mathbb R^{2n+2}_n$,
\begin{equation*}
\begin{array}{c}
J\dfrac{\partial}{\partial x^{\alpha}}=\dfrac{\partial}{\partial x^{\alpha+n+1}},\quad J\dfrac{\partial}{\partial x^{\alpha+n+1}}=\dfrac{\partial}{\partial x^{\alpha}},\\[0.5cm]
G\Big(\dfrac{\partial}{\partial x^{\alpha}},\dfrac{\partial}{\partial x^{\alpha}}\Big)=1,\quad 
G\Big(\dfrac{\partial}{\partial x^{\alpha+n+1}},\dfrac{\partial}{\partial x^{\alpha+n+1}}\Big)=-1
\end{array}
\end{equation*} 
for $\alpha=1,\ldots,n$, where $(x^1,\ldots,x^{2n+2})$ are the Cartesian coordinates on $\mathbb R^{2n+2}_n$.

Consider the hypersurface $H^{2n+1}_n$ in $\mathbb R^{2n+2}_{n+1}$ given by the equation
$$
H^{2n+1}_n=\Big\{x\in R^{2n+2}_{n+1}:\sum_{\alpha=1}^{n+1}x_{\alpha}^2-\sum_{\alpha=n+2}^{2n+2}x_{\alpha}^2=-1\Big\}.
$$
Let 
$$
  N = \sum_{\alpha=1}^{2n+2} x^i \dfrac{\partial}{\partial x^i}
$$
be the normal vector field of $H^{2n+1}_n$. Then $G(N,N)=-1$. Define a vector field $\xi$, a tensor field $\varphi$ of type $(1,1)$, a $1-$form $\eta$ and a pseudo-Riemannian metric $g$ on $H^{2n+1}_n$ by assuming 
\begin{equation*}
  \xi = -JN, \quad JX= \varphi X-\eta(X)\nu,\quad g=G|H^{2n+1}_n.
\end{equation*}
We get an almost paracontact metric structure $(\varphi,\xi,\eta,g)$ on  $H^{2n+1}_n$. We shall show that this structure is para-Sasakian and $H^{2n+1}_n$ is of constant sectional curvature (-1).

For our hypersurface, the Weingarten formula is $D_XN=X$, $D$ being the Levi-Civita connection of $G$. Hence, we obtain the shape operator $A=-I$ and the second fundamental form $h(X,Y)=g(X,Y)$. Using the parallelity of $J$ and the Gauss formula, we have
\begin{eqnarray*}
  0=(D_XJ)\xi=D_XJ\xi-J(D_X\xi) 
     &=& -D_XN-J(\nabla_X\xi+h(X,\xi)N)\\
     &=&-X-\varphi\nabla_X\xi+\eta(X)\xi,
\end{eqnarray*}
where $\nabla$ is the Levi-Civita connection of $g$. Applying $\varphi$ we get $\nabla_X\xi=-\varphi X$. This yields
$$
d\eta(X,Y)=\frac{1}{2}(g(\nabla_X\xi,Y)-g(\nabla_Y\xi,X))=g(X,\varphi Y)
$$
and the structure $(\varphi,\xi,\eta,g)$ is paracontact metric. Moreover, 
\begin{eqnarray*}
  0 &=& (D_XJ)Y=D_X(\varphi Y-\eta(Y)N)-J(\nabla_XY+g(X,Y)N) \\
    &=& \nabla_X\varphi Y+g(X,\varphi Y)N-(X\eta(Y)N-\eta(Y)X \\
    & & \null-\varphi\nabla_XY+\eta(\nabla_XY)N+g(X,Y)\xi \\ 
    &=& (\nabla_X\varphi)Y+g(X,Y)\xi-\eta(Y)X-((\nabla_X\eta)Y-g(X,\varphi Y))N.
\end{eqnarray*}
Taking the tangential part we see that $(\nabla_X\varphi)Y=-g(X,Y)\xi+\eta(Y)X$, and in view of (\ref{sas}) the structure is para-Sasakian.

From the Gauss equation we see that $H^{2n+1}_n$ is the hypersurface of constant curvature equal to 
minus one. Indeed, we have
$$
  R(X,Y)Z=h(Y,Z)AX-h(X,Z)AY=-(g(Y,Z)X-g(X,Z)Y).
$$
\space\hfill$\Box$
\end{ex}

In \cite[Th.\ 3.10]{SZ}, S. Zamkovoy established a certain relationship for the values of the Ricci curvature of a conformally flat para-Sasakian manifold of dimension $2n+1\geqslant 5$. We shall generalize his result and prove the following theorem.

\begin{tw}
A conformally flat para-Sasakian manifold is of constant sectional curvature $k=-1$.
\end{tw}

\begin{pr}
Denote by $\mathop{Ric}$ and $\mathop{Ric^*}$ the Ricci and $*$-Ricci operators, and by $\mathop{r}$ and $\mathop{r}^*$ the scalar and $*$-scalar curvatures, 
\begin{eqnarray*}
  g(\mathop{Ric}Y,Z)&=&\mathop{Tr}\left\{X\rightarrow R(X,Y)Z\right\},\\
  g(\mathop{Ric^*}Y,Z)&=&\mathop{Tr}\left\{X\rightarrow -\varphi R(X,Y)\varphi Z\right\},\\
  \mathop{r}&=&\mathop{Tr} \mathop{Ric},\\
  \mathop{r^*}&=&\mathop{Tr} \mathop{Ric^*}.
\end{eqnarray*}
For a para-Sasakian manifold, the following identities are known (cf. \cite{SZ})
\begin{eqnarray}
\label{Rsasaki}
  & R(X,Y)\xi=\eta(X)Y-\eta(Y)X, & \\
\label{ricsasaki}
  & g(\mathop{Ric}\xi,Y)=-2n\eta(Y), & \\
\label{rr*}
  & \mathop{r}\null + \mathop{r^*}\null + 4n^2 = 0. &
\end{eqnarray}

At first, consider the case $n>1$ ($\dim M>3$). The  Riemann curvature tensor field of a conformally flat manifold is given as  
\begin{eqnarray}
\label{flat}
  R(X,Y)Z&=&\frac{1}{2n-1}(g(Y,Z)\mathop{Ric}X+g(\mathop{Ric}Y,Z)X\nonumber\\
         & &\null-g(X,Z)\mathop{Ric}Y-g(\mathop{Ric}X,Z)Y) \nonumber\\
         & &\null-\frac{r}{2n(2n-1)}(g(Y,Z)X-g(X,Z)Y). 
\end{eqnarray}
Using (\ref{flat}) and (\ref{ricsasaki}), we get
\begin{eqnarray}
\label{r*}
  \mathop{r^{*}} &=& -\frac{1}{2n-1}(\mathop{r}+4n).
\end{eqnarray}
Now, from (\ref{r*}) and (\ref{rr*}), we obtain 
\begin{equation}
\label{r}
  r=-2n(2n+1).
\end{equation}
Taking $Y=Z=\xi$ in (\ref{flat})  and using (\ref{Rsasaki}), (\ref{ricsasaki}) and (\ref{r}), we have $\mathop{Ric}X=-2nX$. Putting the last relation and (\ref{r}) into (\ref{flat}) gives
$$
R(X,Y)Z=-(g(Y,Z)X-g(X,Z)Y), 
$$
which completes the proof in this case. 

Now, consider the case $n=1$ ($\dim M=3$). Since $M$ is conformally flat, we have 
\begin{eqnarray}
\label{flat3}
&&g((\nabla_X\mathop{Ric})Y,Z)-g((\nabla_Z\mathop{Ric})Y,X)\nonumber\\
&&\qquad\qquad\null-\frac{1}{4}((\nabla_Xr)g(Y,Z)-(\nabla_Zr)g(Y,X))=0.
\end{eqnarray}
Since $\xi$ is a Killing vector field, we have
\begin{eqnarray}
\label{lric}
  0 &=& g((\mathcal{L}_{\xi}\mathop{Ric}) Y,Z) \\
    &=& g((\nabla_{\xi}\mathop{Ric})Y,Z)+g(\mathop{Ric}\nabla_{Y}\xi,Z)
        +g(\mathop{Ric}Y,\nabla_Z\xi),\nonumber\\
\label{lr}
  \mathcal{L}_{\xi}r&=&\nabla_{\xi}r=0.
\end{eqnarray}
Using (\ref{flat3}), (\ref{lr}) and (\ref{pScons}), we find 
\begin{eqnarray}
\label{1}
  \quad g((\nabla_{\xi}\mathop{Ric})Y,Z)
  &=&g((\nabla_Z\mathop{Ric})\xi,Y)-\frac{1}{4}(\nabla_Zr)\eta(Y)\\
  &=&\null-2ng(\varphi Y,Z)+g(\mathop{Ric}\varphi Z,Y)-\frac{1}{4}(\nabla_Zr)\eta(Y).\nonumber
\end{eqnarray}
On the other hand, by (\ref{lric}) and (\ref{pScons}), we get  
\begin{eqnarray}
\label{2}
g((\nabla_{\xi}\mathop{Ric})Y,Z)&=&g(\mathop{Ric}\varphi Y,Z)+g(\mathop{Ric}Y,\varphi Z).
\end{eqnarray}
Comparing (\ref{1}) and (\ref{2}) gives
\begin{equation}
g(\mathop{Ric}\varphi Y, Z)=-2ng(\varphi Y,Z)-\frac{1}{4}\eta(Y)\nabla_Zr.
\end{equation}
Puting $\varphi Y$ instead of $Y$ in the above equality and using (\ref{ricsasaki}), we obtain
\begin{equation*}
g(Ric Y,Z)=-2ng(Y,Z).
\end{equation*} 
So $M$ is an Einstein space, and consequently, it is of constant sectional curvature. In view of (\ref{Rsasaki}), the sectional curvature $k=-1$, which completes the proof. 
\hfill $\Box$ 
\end{pr}

To finish this section, we give an example of a class of new paracontact metric structures which are para-CR and not normal in general.

\begin{ex}
\label{p1}
Let $(x^{\alpha},y_{\alpha},z)$, $n\geqslant2$, $1\leqslant\alpha\leqslant n$, be the Cartesian coordinates in $\mathbb R^{2n+1}$. Define a frame of vector fields $(e_i;1\leqslant i\leqslant2n+1)$ on $\mathbb R^{2n+1}$ by 
$$
  e_{\alpha}=\frac{\partial}{\partial x^{\alpha}},\quad
  e_{n+\alpha}=\null-f\frac{\partial}{\partial x^{\alpha}}
               +\frac{\partial}{\partial y_{\alpha}}
               -2x^{\alpha}\frac{\partial}{\partial z},\quad
  e_{2n+1}=\frac{\partial}{\partial z},
$$
$f$ being an arbitrary function. The dual frame of 1-forms $(\theta^i)$ is given by
$$
  \theta^{\alpha}= dx^{\alpha}+f dy_{\alpha},\quad
  \theta^{n+\alpha}= dy_{\alpha},\quad
  \theta^{2n+1}=\null2\sum_{\omega=1}^n x^{\omega}dy_{\omega}+dz.
$$
Define an almost paracontact structure $(\varphi,\xi,\eta,g)$ on $\mathbb R^{2n+1}$ by assuming
\begin{eqnarray*}
  & \varphi e_{\alpha}=-e_{\alpha},\quad 
    \varphi e_{n+\alpha}=e_{n+\alpha},\quad \varphi e_{2n+1}=0,\quad
    \xi=e_{2n+1},\quad \eta=\theta^{2n+1}, &\\
  & g(e_{\alpha},e_{n+\alpha})=g(e_{n+\alpha},e_{\alpha})=g(e_{2n+1},e_{2n+1})=1,\quad\\
    & g(e_i,e_j)=0\quad\mbox{otherwise}. &
\end{eqnarray*}
Here, we have $\varPhi(e_{\alpha},e_{n+\alpha})= -\varPhi(e_{n+\alpha},e_{\alpha})=1$ and $\varPhi(e_{i},e_{j})=0$ otherwise. Therefore the fundamental 2-form $\varPhi$ has the shape 
$$
  \varPhi=
  \sum_{\alpha=1}^n(\theta^{\alpha}\otimes\theta^{n+\alpha}-\theta^{n+\alpha}\otimes\theta^{\alpha})=
  2\sum_{\alpha=1}^n \theta^{\alpha}\wedge\theta^{n+\alpha}=
  2\sum_{\alpha=1}^n dx^{\alpha}\wedge dy_{\alpha}=d\eta.
$$
Thus, $(\varphi,\xi,\eta,g)$ is a paracontact metric structure. Let us calculate the structure tensor $h$.
\begin{eqnarray*}
he_{\alpha}&=&\frac{1}{2}\Big(\Big[\xi,\varphi e_{\alpha}\Big]-\varphi\Big[\xi, e_{\alpha}\Big]\Big)=\frac{1}{2}\Big(-\Big[\frac{\partial}{\partial z},\frac{\partial}{\partial x^{\alpha}}\Big]-\varphi\Big[\frac{\partial}{\partial z},\frac{\partial}{\partial x^{\alpha}}\Big]\Big)=0;\\
he_{n+\alpha}&=&\cfrac{1}{2}\Big(\Big[\xi,\varphi e_{n+\alpha}\Big]-\varphi\Big[\xi, e_{n+\alpha}\Big]\Big)\\
&=&\cfrac{1}{2}\Big(\Big[\cfrac{\partial}{\partial z},-f\cfrac{\partial}{\partial x^{\alpha}}+\cfrac{\partial}{\partial y_{\alpha}}-2x^{\alpha}\cfrac{\partial}{\partial z}\Big]
-\varphi \Big[\cfrac{\partial}{\partial z},-f\cfrac{\partial}{\partial x^{\alpha}}
+\cfrac{\partial}{\partial y_{\alpha}}-2x^{\alpha}\cfrac{\partial}{\partial z}\Big]\Big)\\
&=&\cfrac{1}{2}\Big(-\cfrac{\partial f}{\partial z}\cfrac{\partial}{\partial x^{\alpha}}+\varphi\cfrac{\partial f}{\partial z}\cfrac{\partial}{\partial x^{\alpha}}\Big)
=-\cfrac{\partial f}{\partial z}e_{\alpha};\\
he_{2n+1}&=&h\xi=0.
\end{eqnarray*}
We see that $h^2=0$, and $h\neq0\Leftrightarrow \frac{\partial f}{\partial z}\neq0$. Therefore, if $\frac{\partial f}{\partial z}\neq0$, the structure is not para-Sasakian one.
To find necessary and sufficient conditions for this structure to be para-CR, we are going to use the condition (\ref{plusminus}). We see that the eigendistribution $\mathcal D^{-}$ is spanned by the vector fields $(e_{\alpha})$. Since $[e_{\alpha},e_{\beta}]=0$, $\mathcal D^{-}$ is involutive. The eigendistribution $\mathcal D^{+}$ is spanned by the vector fields $(e_{n+\alpha})$. Since
$$
  [e_{n+\alpha},e_{n+\beta}]=
  \Big(f\frac{\partial f}{\partial x^{\alpha}}
  -\frac{\partial f}{\partial y_{\alpha}}
  +2x^{\alpha}\frac{\partial f}{\partial z}\Big)e_{\beta}
  -\Big(f\frac{\partial f}{\partial x^{\beta}}
  -\frac{\partial f}{\partial y_{\beta}}
  +2x^{\beta}\frac{\partial f}{\partial z}\Big)e_{\alpha},
$$
$\mathcal D^{+}$ is involutive if and only if the function $f$ fulfills the following system of partial differential equations
\begin{equation}
\label{nsex1}
  f\frac{\partial f}{\partial x^{\alpha}}
  -\frac{\partial f}{\partial y_{\alpha}}
  +2x^{\alpha}\frac{\partial f}{\partial z}=0,
  \quad 1\leqslant\alpha\leqslant n.
\end{equation}
Therefore, $(\varphi,\xi,\eta,g)$ is a para-CR structure if and only if (\ref{nsex1}) holds. 

To see that this structure is not normal in general, at first we find 
$$
  [e_{n+\alpha},\xi]=\frac{\partial f}{\partial z} e_{\alpha},
$$
and next we compute
\begin{eqnarray*}
  N(e_{n+\alpha},\xi)-2d\eta(e_{n+\alpha},\xi)\xi
    &=& \varphi^2[e_{n+\alpha},\xi]-\varphi[\varphi e_{n+\alpha},\xi]-2d\eta(e_{n+\alpha},\xi)\xi \\
    &=& [e_{n+\alpha},\xi]-\eta([e_{n+\alpha},\xi])\xi-\varphi[e_{n+\alpha},\xi] \\
    &=& 2\frac{\partial f}{\partial z} e_{\alpha},
\end{eqnarray*}
which is non-zero if $\frac{\partial f}{\partial z}\not=0$. 

Finally, note that the function 
$$
  f(x^{\alpha},y_{\alpha},z)=\frac1z\Big(c+\sum_{\omega=1}^n (x^{\omega})^2\Big),\quad
  c=\mbox{const.},
$$
is an example of a particular solution of (\ref{nsex1}), and moreover $\frac{\partial f}{\partial z}\not=0$. 
\hfill $\Box$
\end{ex}

\section*{ \normalsize 5.  Almost Para-Cosymplectic Manifolds}

\noindent
Before we study almost paracosymplectic manifolds, recall the notion of para-K\"ahlerian manifolds and almost para-K\"ahlerian manifolds; cf. \cite{CFG,CGM}.

An almost para-K\"ahlerian manifold $\widetilde{M}$ by definition is a $2n$-dimensional differentiable manifold endowed with an almost para-K\"ahlerian structure $(\widetilde{J},\widetilde{g})$. The structure is formed by a (1,1)-tensor field $\widetilde{J}$ such that $\widetilde{J}^2=\widetilde{I}$, and a pseudo-Riemannian metric $\widetilde{g}$ satysfying $\widetilde{g}(\widetilde{J}\widetilde{X},\widetilde{J}\widetilde{Y})=-\widetilde{g}(\widetilde{X},\widetilde{Y})$, and the fundamental form $\widetilde\varPhi$, $\widetilde{\varPhi}(\widetilde{X},\widetilde{Y}) = \widetilde{g}(\widetilde{X},\widetilde{J}\widetilde{Y})$, is closed. An almost para-K\"ahlerian manifold with integrable almost para-complex structure $\widetilde{J}$ (equivalently, $\widetilde{\nabla}\widetilde{J}=0$) is said to be para-K\"ahlerian.

Let $M$ be a paracosymplectic manifold. Since $d\eta=0$, the paracontact distribution $\mathcal{D}$ is completely integrable and the manifold $M$ possesses a foliation $\mathcal{F}$ generated by  $\mathcal{D}$. Any leaf $\widetilde{M}$ of $\mathcal{F}$ is a submanifold of $M$ of codimension $1$. Since $\xi|_{\widetilde{M}}$ is a vector field normal to $\widetilde{M}$, we may treat $\widetilde{M}$ as a pseudo-Riemannian hypersurface. Let $\widetilde{J}$ be the (1,1)-tensor field defined by $\widetilde{J}=\varphi|_{\widetilde{M}}$ and $\widetilde{g}$ the induced metric on $\widetilde{M}$. Then the pair $(\widetilde{J},\widetilde{g})$ is an almost para-K\"ahlerian structure on $\widetilde{M}$ (its fundamental form is closed since it is the pullback of the fundamental form of $M$). We say that $M$ has para-K\"ahlerian leaves if on every leaf $\widetilde M$, the induced structures $(\widetilde{J},\widetilde{g})$ are para-K\"ahlerian.

In the paper (\cite{D}), P. Dacko shows that an almost paracoymplectic manifold satisfies the following conditions  
\begin{eqnarray}
\label{dackowzorki}
  & \nabla_{\xi}\xi=0, \quad \nabla_{\xi}\varphi=0, \quad \nabla_{\varphi X}\xi=-\varphi\nabla_X\xi, &\\
\label{dacko1}
  & (\nabla_{\varphi X}\varphi)\varphi Y-(\nabla_X\varphi)Y=\eta(Y)\varphi\nabla_{ X}\xi. &
\end{eqnarray}
Morover, he proves that an almost para-cosympletic manifold  has para-K\"ahlerian leaves if and only if 
\begin{equation}
\label{wzor2}
  (\nabla_X\varphi)Y=g(\varphi\nabla_X\xi,Y)\xi-\eta(Y)\varphi\nabla_X\xi.
\end{equation}
Using the above equality, we obtain  

\begin{tw}
An almost para-cosympletic manifold is a para-CR manifold if and only if it has para-K\"ahlerian leaves.
\end{tw}

\begin{pr}
We prove that for a para-cosymplectic manifold, formulas (\ref{wzor2}) and (\ref{paracr}) are equivalent to each other. 

To do it, we check at first that under (\ref{dackowzorki}), (\ref{dacko1}) and the general formula 
$\allowbreak(\nabla_U\eta)(V) =(\nabla_V\eta)(U)$ (which is a consequence of $d\eta=0$), formula (\ref{paracr}) reduces equivalently to 
\begin{equation}
\label{paracrcos}
 (\nabla_X\varphi)Y = g(\varphi\nabla_X\xi,Y)\xi
\end{equation}
for any $X,Y\in \mathcal{D}$. We see that (\ref{wzor2}) implies (\ref{paracrcos}). Conversely, let us assume that (\ref{paracrcos}) is fulfilled. Let $X,Y\in\mathfrak{X}(M)$. Then $X-\eta(X)\xi\in\mathcal D$ and $Y-\eta(Y)\xi\in\mathcal D$. Replacing $X$ by $X-\eta(X)\xi$ and $Y$ by $Y-\eta(Y)\xi$ in (\ref{paracrcos}), and applying (\ref{dackowzorki}), we find
$$
  (\nabla_X\varphi)Y - \eta(Y)(\nabla_X\varphi)\xi = g(\varphi\nabla_X\xi,Y)\xi,
$$
which leads to (\ref{wzor2}). 
\hfill $\Box$
\end{pr}

Below, we give a class of new almost paracosymplectic structures which are para-CR and not normal in general.

\begin{ex}
Let $(x^{\alpha},y_{\alpha},z)$ ($n\geqslant1$, $1\leqslant\alpha\leqslant n$) be the Cartesian coordinates in $\mathbb R^{2n+1}$. Let for $1\leqslant\alpha,\beta\leqslant n$, $F^{\alpha}_{\beta}\colon\mathbb R^{2n+1}\to\mathbb R$ be functions depending on the coordinates $x^{\alpha}$ and $z$ only, and such that $F^{\beta}_{\alpha}=F^{\alpha}_{\beta}$. Define a frame of vector fields $(e_i;1\leqslant i\leqslant2n+1)$ in $\mathbb R^{2n+1}$ by 
$$
  e_{\alpha}=\frac{\partial}{\partial x^{\alpha}}
        -\sum_{\omega=1}^n F^{\omega}_{\alpha}\frac{\partial}{\partial y_{\omega}},\quad
  e_{n+\alpha}=\frac{\partial}{\partial y_{\alpha}},\quad
  e_{2n+1}=\frac{\partial}{\partial z}.
$$
Then the dual frame of 1-forms $(\theta^i;1\leqslant i\leqslant2n+1)$ is given by
$$
  \theta^{\alpha}= dx^{\alpha},\quad
  \theta^{n+\alpha}= dy_{\alpha} + \sum_{\omega=1}^n F^{\alpha}_{\omega}dx^{\omega},\quad
  \theta^{2n+1}=dz.
$$
Define an almost paracontact structure $(\varphi,\xi,\eta,g)$ on $\mathbb R^{2n+1}$ by assuming
\begin{eqnarray*}
  & \varphi e_{\alpha}=-e_{\alpha},\quad 
    \varphi e_{n+\alpha}=e_{n+\alpha},\quad 
    \varphi e_{2n+1}=0,\quad
    \xi=e_{2n+1},\quad 
    \eta=\theta^{2n+1}, \\
  & g(e_{\alpha},e_{n+\alpha})=g(e_{n+\alpha},e_{\alpha})=g(e_{2n+1},e_{2n+1})=1. &
\end{eqnarray*}
Here, we have $\varPhi(e_{\alpha},e_{n+\alpha})= -\varPhi(e_{n+\alpha},e_{\alpha})=1$ and $\varPhi(e_{i},e_{j})=0$ otherwise. Therefore the fundamental 2-form $\varPhi$ has the shape 
$$
  \varPhi = 2\sum_{\alpha=1}^n \theta^{\alpha}\wedge\theta^{n+\alpha} 
          = 2\sum_{\alpha=1}^n dx^{\alpha}\wedge dy_{\alpha}
             +2\sum_{\alpha,\omega=1}^n F^{\alpha}_{\omega} dx^{\alpha}\wedge dx_{\omega}.
$$
By the symmetry $F^{\alpha}_{\omega}=F_{\alpha}^{\omega}$, it must be that 
$$
  \varPhi = 2\sum_{\alpha=1}^n dx^{\alpha}\wedge dy_{\alpha},
$$
and consequently, $d\varPhi=0$. Since also $d\eta=0$, the structure $(\varphi,\xi,\eta,g)$ is almost paracosymplectic.

As in the previous example, to find necessary and sufficient conditions for this structure to be para-CR, we shall use the condition (\ref{plusminus}). The eigendistribution $\mathcal D^{+}$ is spanned by the vector fields $(e_{n+\alpha})$. Since $[e_{n+\alpha},e_{n+\beta}]=0$, $\mathcal D^{+}$ is involutive. The eigendistribution $\mathcal D^{-}$ is spanned by the vector fields $(e_{\alpha})$. Since
$$
  [e_{\alpha},e_{\beta}] = \sum_{\omega}^n 
    \Bigg(\frac{\partial F^{\omega}_{\alpha}}{\partial x^{\beta}}
    -\frac{\partial F^{\omega}_{\beta}}{\partial x^{\alpha}}\Bigg) \frac{\partial}{\partial y^{\omega}},
$$
$\mathcal D^{-}$ is involutive if and only if 
$$
  \frac{\partial F^{\omega}_{\alpha}}{\partial x^{\beta}}
    -\frac{\partial F^{\omega}_{\beta}}{\partial x^{\alpha}} = 0,
$$
or equivalently, 
$$
  F_{\alpha}^{\omega} = \frac{\partial G^{\omega}}{\partial x^{\alpha}}
$$
for certain functions $G^{\omega}$ on $\mathbb R^{2n+1}$ depending on $x^{\alpha}$ and $z$ only. By $F_{\alpha}^{\omega}=F^{\alpha}_{\omega}$, 
$$
  \frac{\partial G^{\omega}}{\partial x^{\alpha}}-\frac{\partial G^{\omega}}{\partial x_{\alpha}}=0
$$
and consequently, 
$$
  G^{\omega} = \frac{\partial H}{\partial x^{\omega}}
$$
for a certain function $H$ on $\mathbb R^{2n+1}$ depending on $x^{\alpha}$ and $z$ only. Thus, 
$\mathcal D^{-}$ is involutive if and only if
$$
  F_{\alpha}^{\beta} = \frac{\partial^2 H}{\partial x^{\alpha}\partial x^{\beta}},
$$
$H$ being a function $H$ on $\mathbb R^{2n+1}$ depending on $x^{\alpha}$ and $z$ only.

We check that this structure is not normal. At first, we find 
$$
  [e_{\alpha},\xi]=
  \sum_{\omega=1}^n \frac{\partial^2 H}{\partial x^{\omega}\partial x^{\alpha}{\partial z}} e_{n+\omega},
$$
and next we compute
\begin{eqnarray*}
  N(e_{\alpha},\xi)-2d\eta(e_{\alpha},\xi)\xi
    &=& \varphi^2[e_{\alpha},\xi]-\varphi[\varphi e_{\alpha},\xi]-2d\eta(e_{\alpha},\xi)\xi \\
    &=& [e_{\alpha},\xi]-\eta([e_{\alpha},\xi])\xi+\varphi[e_{\alpha},\xi] \\
    &=& 2\sum_{\omega=1}^n \frac{\partial^2 H} 
                  {\partial x^{\omega}\partial x^{\alpha}{\partial z}} e_{n+\omega},
\end{eqnarray*}
which is non-zero in general. 
\end{ex}

\smallskip
\noindent
\textbf{Note}\\[6pt]
\noindent
Certain of the presented results can be treated as paracontact analogies of those known in contact geometry. For contact metric manifolds, we refer the beautiful book \cite{B}.

\bigskip
\noindent
\textbf{Acknowledgments}\\[6pt]
\noindent
I would like to thank my supervisor prof. Zbigniew Olszak for his support
and suggestions during this work.


\begin{thebibliography}{1}
\small
\bibitem{AMT1}
D. V. Alekseevsky, C. Medori and A. Tomassini,
{\it Maximally homogeneous para-CR manifolds}
Ann. Global Anal. Geom. {\bf 30} (2006), No. 1, 1-27.

\bibitem{AMT2}
D. V. Alekseevsky, C. Medori and A. Tomassini,
{\it Maximally homogeneous para-CR manifolds of semisimple type},
in: {\it Handbook in Pseudo-Riemannian Geometry and Supersymmetry}, 
(IRMA Lectures in Mathematics and Theoretical Physics 16, 
Z\"urich: European Mathematical Society, 2010) 569-577.

\bibitem{B}
D. E. Blair,
{\it Riemannian geometry of contact ansd symplectic manifolds},
Progress in Mathematics Vol. {\bf203}, Birkh\"auser, Boston (Massachusates), 2002.

\bibitem{CM}
B. Cappelletti Montano,
{\it Bi-Legendrian structures and paracontact geometry},
Int. J. Geom. Meth. Mod. Phys. 6 (2009), 487-504.

\bibitem{CFG}
V. Cruceanu, P. Fortuny and P. M. Gadea,
{\it A survey on paracomplex geometry},
Rocky Mountain J. Math. {\bf 26} (1996), 83-115.

\bibitem{CGM}
V. Cruceanu, P. M. Gadea and J. Mu\~noz Masqu\'e,
{\it Para-Hermitian and para-K\"ahler manifolds},
Quaderni dell`Istituto di Matematica, Facolt\`a di Economia, 
Universit\`a di Messina, No. {\bf 1} (1995), 1-72.

\bibitem{D}
P. Dacko,
{\it On almost para-cosymplectic manifolds},
Tsukuba J. Math. {\bf 28} (2004), No. 1, 193-213.

\bibitem{DT}
S. Dragomir and G. Tomassini,
{\it Differential geometry and analysis on CR manifolds},
Progress in Mathematics Vol. {\bf246}, Birkh\"auser, Boston (Massachusates), 2006.

\bibitem{E}
S. Erdem,
{\it On almost (para)contact (hiperbolic) metric manifolds and harmonicity of $(\varphi,\varphi')$-holomorphic maps between them},
Houston J. Math. {\bf 28} (2002), No. 1, 21-45. 

\bibitem{HN}
C. D. Hill and P. Nurowski,
{ Differential equations and para-CR structures},
{\it Boll. Unione Mat. Ital.} (9) {\bf3} (2010), 25-91.

\bibitem{IVZ}
S. Ivanov, D. Vassilev and S. Zamkovoy,
{\it Conformal paracontact curvature and the local flatness theorem},
Geometriae Dedicata, {\bf144} (2010), No. 1, 79-100.

\bibitem{K}
S. Kaneyuki, 
{\it On classification of para-Hermitian symmetric spaces}, 
Tokyo J. Math. {\bf8} (1985), 473-482.

\bibitem{KW}
S. Kaneyuki and F. L. Williams,
{\it Almost paracontact and paraHodge structures on manifolds},
Nagoya Math. J. {\bf 99} (1985), 173-187.

\bibitem{Kush}
A. Kushner,
{ Almost product structures and Monge-Amp\`ere equations},
{\it Lobachevskii J. Math.} {\bf23} (2006), 151-181.

\bibitem{NS}
P. Nurowski and G. A. J. Sparling,
{ Three-dimensional Cauchy–Riemann structures and second-order ordinary differential equations},
{\it Class. Quantum Grav.} {\bf20} (2003), 4995-5016.

\bibitem{JW}
J. We{\l}yczko,
{\it On Legendre curves in 3-dimensional normal almost paracontact metric manifolds},
Results Math. {\bf54} (2009), 377-387.

\bibitem{SZ}
S. Zamkovoy,
{\it Canonical connections on paracontact manifolds}, 
Ann. Glob. Anal. Geom. {\bf36} (2009), No. 1, 37-60.
\end{thebibliography}
\end {document}